\documentclass[11pt]{article}
\usepackage{amssymb,latexsym}
\usepackage{epsfig}
\usepackage{eufrak}
\usepackage{amsmath}
\usepackage{mathrsfs}
\usepackage{color}

\newtheorem{theorem}{Theorem}
\newtheorem{lemma}{Lemma}

\newtheorem{remark}{Remark}

\newcommand{\be}{\begin{equation}}
\newcommand{\ee}{\end{equation}}
\newcommand{\bee}{\begin{eqnarray*}}
\newcommand{\eee}{\end{eqnarray*}}
\newcommand{\bel}{\begin{eqnarray}}
\newcommand{\eel}{\end{eqnarray}}
\newcommand{\bec}{\begin{cases}}
\newcommand{\eec}{\end{cases}}
\newcommand{\bem}{\begin{bmatrix}}
\newcommand{\eem}{\end{bmatrix}}

\newcommand{\li}{\left}
\newcommand{\ri}{\right}

\newcommand{\DEF}{\stackrel{\mathrm{def}}{=}}

\newcommand{\ep}{\epsilon}
\newcommand{\vep}{\varepsilon}

\newcommand{\si}{\sigma}
\newcommand{\Si}{\Sigma}

\newcommand{\ro}{\rho}

\newcommand{\f}{\frac}
\newcommand{\sq}{\sqrt}
\newcommand{\cd}{\cdots}

\newcommand{\qqu}{\qquad}

\newcommand{\mrm}{\mathrm}

\newcommand{\bed}{\begin{description}}
\newcommand{\eed}{\end{description}}
\newcommand{\bei}{\begin{itemize}}
\newcommand{\eei}{\end{itemize}}
\newcommand{\ben}{\begin{enumerate}}
\newcommand{\een}{\end{enumerate}}
\newcommand{\bib}{\bibitem}
\newcommand{\beL}{\begin{lemma}}
\newcommand{\eeL}{\end{lemma}}
\newcommand{\beT}{\begin{theorem}}
\newcommand{\eeT}{\end{theorem}}

\newcommand{\bpf}{\begin{pf}}
\newcommand{\epf}{\end{pf}}
\newcommand{\bsk}{\bigskip}

\newcommand{\pfbox}{\hfill\mbox{$\Box$}}

\newenvironment{pf}{\paragraph*{Proof{\rm.}}}{\pfbox\bigskip}

\begin{document}

\title{{\bf A New Generalization of Chebyshev Inequality for Random Vectors}
\thanks{The author is with
Department of Electrical and Computer Engineering, Louisiana State
University, Baton Rouge, LA 70803; Email: chenxinjia@gmail.com. }}

\author{Xinjia Chen}

\date{June 2007}

\maketitle


\begin{abstract}

In this article, we derive a new generalization of Chebyshev inequality for random vectors. We demonstrate that
the new generalization is much less conservative than the classical generalization.

\end{abstract}

\bigskip

\section{Classical Generalization of Chebyshev inequality}

The Chebyshev inequality discloses the fundamental relationship between the mean and variance of a random
variable. Extensive research works have been devoted to its generalizations for random vectors. For example,
various generalizations can be found in Marshall and Olkin (1960), Godwin (1955), Mallows (1956) and the
references therein. A natural generalization of Chebyshev inequality  is as follows. \bsk

For a random vector $X \in {\bf R}^{n}$ with cumulative distribution $F(.)$, \be \Pr\left\{\left||X-E[X]\right||
\geq \varepsilon\right\} \leq \frac{{\rm Var}(X)}{{\varepsilon}^{2}} \;\;\;\; \forall \varepsilon > 0
\label{eq1} \ee
 where $||.||$ denotes the Euclidean norm of a vector and
\[
{\rm Var} (X) \DEF \int_{ V \in \; {\bf R}^{n} } ||V - E[X]||^2 dF(V)
\]
This classical generalization can be found in a number of textbooks of probability theory and statistics (see,
e.g., pp. 446-451 of Laha and Rohatgi (1979)).

\section{New Generalization of Chebyshev inequality}

The classical generalization (\ref{eq1}) perfectly assembles its counterpart for scalar random variables.
However, it may be too conservative.  To address the conservatism,  we derive a new multivariate Chebyshev
inequality as follows.

\begin{theorem} \label{gc}
For any random vector $X \in {\bf R}^{n}$ with covariance matrix $\Si$, \be \Pr\left\{ (X-E[X])^\top
\Sigma^{-1}(X-E[X]) \geq \varepsilon \right\} \leq \frac{n}{\varepsilon}, \qqu \forall{\varepsilon}
> 0 \label{eq2} \ee where the superscript ``$\top$'' denotes the transpose of a matrix.
\end{theorem}

\begin{pf}
Let ${\bf D}_{\varepsilon}=\left\{V \in {\bf R}^{n}:\;(V-E[X])^\top \Sigma^{-1}(V-E[X]) \geq \varepsilon
\right\}$. By the definition of ${\bf D}_{\varepsilon}$, we have
\[
\frac{1}{\varepsilon}(V-E[X])^\top \Sigma^{-1}(V-E[X]) \geq 1, \qqu \forall V \in {\bf D}_{\varepsilon}.
\]
Hence,
\begin{eqnarray*}
\Pr\left\{X\in{{\bf D}}_{\varepsilon}\right\} & \leq & \frac{1}{\varepsilon} \; {\int}_{ V \in \; {\bf
D}_{\varepsilon}} (V-E[X])^\top
\Sigma^{-1}(V-E[X]) dF(V)\\
& \leq &  \frac{1}{\varepsilon} \; {\int}_{ V \in \; {\bf R}^{n} } (V-E[X])^\top \Sigma^{-1} (V-E[X]) dF(V).
\end{eqnarray*}
For $i = 1, \cd, n$, let $u_i$ denote the $i$-th element of $V-E[X]$.  For $i=1, \cd, n$ and $j = 1,\cd, n$, let
$\si_{ij}$ denote the element of $\Si$ in the $i$-th row and $j$-th column.  Similarly, let $\ro_{ij}$ denote
the element of $\Si^{-1}$ in the $i$-th row and $j$-th column. Then, \bee (V-E[X])^\top \Sigma^{-1} (V-E[X])
 & = & \sum_{i=1}^n u_i \li ( \sum_{k=1}^n \ro_{ik} \; u_k \ri )\\
& = & \sum_{i=1}^n \sum_{k=1}^n \ro_{ik} \; u_i \; u_k. \eee It follows that \bee &  & {\int}_{ V \in \; {\bf
R}^{n} } (V-E[X])^\top \Sigma^{-1} (V-E[X]) \; dF(V)\\
 & = & {\int}_{ V \in \; {\bf R}^{n} } \li ( \sum_{i=1}^n
\sum_{k=1}^n \ro_{ik} \; u_i \; u_k \ri ) dF(V)\\
& = & \sum_{i=1}^n \sum_{k=1}^n \ro_{ik} \li [ {\int}_{ V \in \; {\bf R}^{n} }  u_i \; u_k  \; dF(V) \ri ].
 \eee
By the definition of the covariance matrix $\Si$ and its symmetry, we have
\[
{\int}_{ V \in \; {\bf R}^{n} }  u_i \; u_k  \; dF(V) = \si_{ik} = \si_{ki}
\]
for $i = 1, \cd, n$ and $k = 1, \cd, n$. Hence, \bee &   & {\int}_{ V \in \; {\bf R}^{n} } (V-E[X])^\top \Sigma^{-1} (V-E[X]) dF(V)\\
 & = & \sum_{i=1}^n
\sum_{k=1}^n \ro_{ik} \si_{ki}\\
& = & \mrm{tr} (\Si^{-1} \Si)\\
& = & n \eee where $\mrm{tr}(.)$ denotes the trace of a matrix.   Therefore,
\begin{eqnarray*}
\Pr\left\{X\in{{\bf D}}_{\varepsilon}\right\}
 & \geq & \frac{1}{\varepsilon} \; {\int}_{ V \in \; {\bf R}^{n} }
 \left(\Sigma^{-1}(V-E[X])(V-E[X])^\top \right) dF(V) \\
& = & \frac{n}{\varepsilon}.
\end{eqnarray*}
The proof is thus completed.
\end{pf}

\begin{remark}
Theorem~\ref{gc} indicates a fundamental relationship between the mean and covariance of a random vector and
describes how a random vector deviates from its expectation.  Specially, for $n = 1$, we have $\Si =
\mrm{Var}(X)$ and by Theorem~\ref{gc}, for any $\ep > 0$, \bee &   & \Pr\left\{ (X-E[X])^\top \Sigma^{-1}(X-E[X]) > \ep \right\}\\
& = & \Pr \li \{ ||X-E[X]|| > \sq{\ep \; \mrm{Var}(X)} \ri \}\\
& \leq & \f{1}{\ep},   \eee from which we deduce
\[
\Pr\left\{\left||X-E[X]\right|| > \varepsilon\right\} \leq \frac{{\rm Var}(X)}{{\varepsilon}^{2}}
\]
by letting $\vep = \sq{\ep \; \mrm{Var}(X)}$. This shows that Theorem~\ref{gc} includes the well-known {\it
Chebyshev} inequality as a special case.
\end{remark}

\begin{remark}

Actually, we had established Theorem~\ref{gc} in \cite[pp.
8--9]{Chen1} in 1997. The applications of this result to control
engineering can be found in \cite{Chen1, Chen2}.  Recently,
Theorem~\ref{gc} has been extended to random elements taking values
in a separate Hilbert space by Rao \cite{Rao} and to random elements
taking values in a separate Banach space by Zhou and Hu \cite{Zhou}.

\end{remark}

\section{Comparison with Classical Generalization}

In this section, we shall show that the inequality in Theorem~\ref{gc} can be much less conservative than the
classical generalized Chebyshev inequality~(\ref{eq1}).

\bsk

 Let $\delta \in (0,1)$.  Based on inequality~(\ref{eq1}), sphere
\[
{\bf B}_{\delta} \DEF \left\{ V \in {\bf R}^n: \left||V-E[X]\right||^2 \leq \frac{ {\rm tr} (\Sigma) } {\delta}
\right\}
\]
is the smallest set
that can be constructed to ensure $\Pr \{ X \in {\bf B}_{\delta}  \} > 1 - \delta$.
On the other hand, by applying Theorem~\ref{gc} we can construct an ellipsoid
\[
{\bf E}_{\delta} \DEF \left\{ V \in {\bf R}^n : (V-E[X])^\top \Sigma^{-1}(V-E[X])  \leq \frac{n}{\delta}
\right\},
\]
 which guarantees $ \Pr \{ X \in {\bf E}_{\delta}  \} > 1 - \delta$.

 For a comparison of the conservativeness of
generalized Chebyshev inequalities~(\ref{eq1}) and~(\ref{eq2}), it is natural to consider the ratio $\frac{ {\rm
vol} ({\bf B}_{\delta} )} { {\rm vol} ({\bf E}_{\delta} )}$ where ${\rm vol} (.)$ is a volume function such that
${\rm vol}(S) = \int_{v \in S} dv$ for any $S \subset {\bf R}^n$.  Interestingly, we have

\begin{theorem} \label{gcc} For any random vector $X \in {\bf R}^n$,
\[
\frac{ {\rm vol} ({\bf B}_{\delta} )} { {\rm vol} ({\bf E}_{\delta} )} =
\frac{ \left( \sqrt{ \frac{ {\rm tr} (\Sigma)   } { n } } \right)^n} { \sqrt{\det(\Sigma)  }  }> 1
\]
where $\det(\Sigma)$ is the determinant of $\Sigma$.
\end{theorem}

\begin{pf}
By the definitions of  variance and covariance,  we have ${\rm Var} (X) = {\rm tr} (\Sigma)$. It follows that
\[
{\rm vol} ({\bf B}_{\delta} ) = K \left( \sqrt{ \frac{ {\rm tr} (\Sigma) } {\delta}  } \right)^n
\]
where $K > 0$ is a constant.  Applying a linear transform $u = \Sigma^{-\frac{1}{2}}(v-E[X]) $ to the
integration ${\rm vol} ({\bf E}_{\delta} ) = \int_{v \in {\bf E}_{\delta}} dv$, we have
\[
{\rm vol} ({\bf E}_{\delta} ) = \det( \Sigma^{\frac{1}{2}}  ) \; \int_{||u||^2 \leq \frac{n}{\delta}} du =
\sqrt{\det({\Sigma})} \; K \left( \sqrt{  \frac{ n } {\delta} } \right)^n
\]
and thus
\[
\frac{ {\rm vol} ({\bf B}_{\delta} )} { {\rm vol} ({\bf E}_{\delta} )} = \frac{ \left( \sqrt{ \frac{ {\rm tr}
(\Sigma)   } { n } } \right)^n} { \sqrt{\det(\Sigma)  }  }. \]
 To show $\frac{ {\rm vol} ({\bf B}_{\delta} )} {
{\rm vol} ({\bf E}_{\delta} )} > 1$, it is equivalent to show
\[
\frac{{\rm tr} (\Sigma)  } {n}  \geq [\det(\Sigma)]^{\frac{1}{n}}. \]
 Recall that the geometric average is no less
than the arithmetic average, \be \frac{{\rm tr} (\Sigma)  } {n} = \frac{\sum_{i=1}^n \sigma_{ii}} {n} \geq
\left( \prod_{i=1}^n \sigma_{ii} \right)^{\frac{1}{n}}, \label{eq3} \ee where $\sigma_{ii}, \;i = 1, \cdots,n$
are the diagonal components of $\Sigma$. Note that the covariance matrix $\Sigma$ is positive definite, hence by
Hadamard's inequality, \be \det(\Sigma) \leq \prod_{i=1}^n \sigma_{ii}. \label{eq4} \ee It follows from
(\ref{eq3}) and (\ref{eq4}) that $ \frac{{\rm tr} (\Sigma)  } {n}  \geq [\det(\Sigma)]^{\frac{1}{n}}$.  The
proof is thus completed.
\end{pf}

As an illustrative example, consider a two-dimensional random vector \[ X = \left[ \begin{array}{l}
y\\
y+z \end{array}  \right] \]
 where $y$ and $z$ are independent Guassian random variables with zero means and
variances $\sigma^2, \; k \sigma^2$ respectively. Straightforward computation gives
\[
\Sigma = \left[
\begin{array} {ll}
\sigma^2 &  \sigma^2\\
\sigma^2 & (k +1)\sigma^2 \end{array} \right] \] and
\[
\frac{ {\rm vol} ({\bf B}_{\delta} )} { {\rm vol} ({\bf
E}_{\delta} )} = \frac{k + 2} {2 \sqrt{k}} \geq \sqrt{2}.
\]
  Obviously, as $k$ increases from $2$ to $\infty$ or decreases from $2$ to $0$, the ratio of
 volumes increases monotonically and tends to $\infty$.

 In the following Figure~\ref{fig}, ellipsoid ${\bf E}_{\delta}$ and
sphere ${\bf B}_{\delta}$ are constructed for $\sigma = 1, \; k = 25$ and $\delta = 0.1$.  Moreover, $1000$
i.i.d. samples of $X$ are generated to show the coverage of the ellipsoid and sphere. It can be seen that most
samples are included in the ellipsoid.  This indicates that Theorem~\ref{gc} is much less conservative than the
classical generalized Chebyshev inequality in describing how a random vector deviates from its expectation.

\begin{figure}[htbp]
\centerline{\psfig{figure=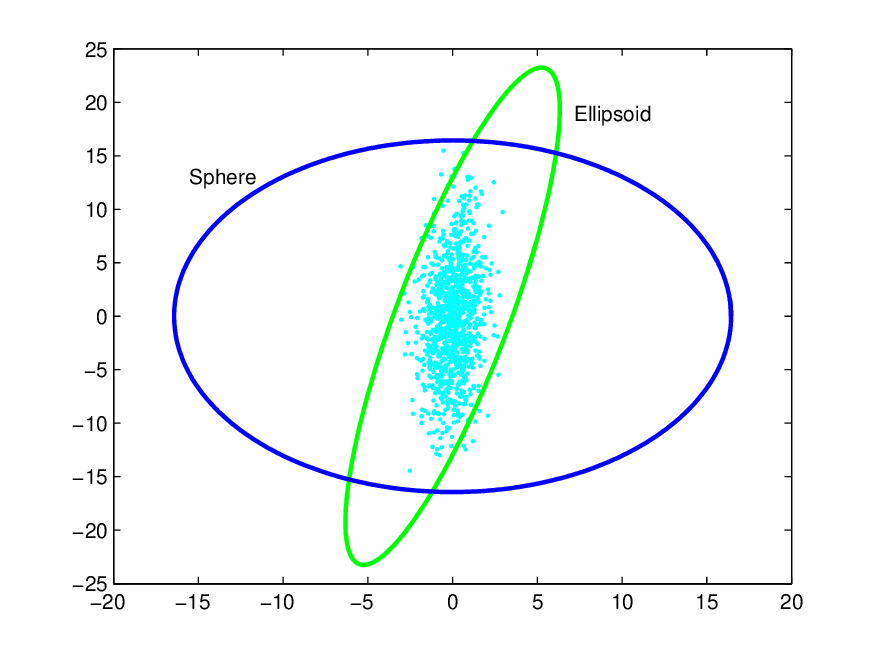, height= 4 in, width= 4 in }} \caption{ Comparison of Generalized
Chebyshev Inequalities } \label{fig}
\end{figure}

\end{document}